\documentclass{article}
\usepackage{amsfonts,amsmath,amssymb,theorem,epsfig,array}

\font\teneusm=eusm10
\font\seveneusm=eusm7
\font\fiveeusm=eusm5
\newfam\eusmfam
\textfont\eusmfam\teneusm
\scriptfont\eusmfam\seveneusm
\scriptscriptfont\eusmfam\fiveeusm
\def\scr#1{{\fam\eusmfam\relax#1}}

\newtheorem{Thm}{Theorem}

\newtheorem{Prop}[Thm]{Proposition}
\newtheorem{Lemma}[Thm]{Lemma}
\newtheorem{Lemma-Def}[Thm]{Lemma-Definition}

{\theorembodyfont{\rmfamily} \newtheorem{Rem}[Thm]{Remark} }
{\theorembodyfont{\rmfamily} \newtheorem{Ex}[Thm]{Example} }
{\theorembodyfont{\rmfamily} \newtheorem{Def}[Thm]{Definition} }

\def\pf{\noindent {\it Proof.}\hskip 8pt}
\def\qed{\hfill{\setlength{\fboxrule}{0.8pt}\setlength{\fboxsep}{1mm}
\fbox{\null}} \vskip 10pt}

\def\frak{\mathfrak}

\newcommand{\sdir}{\leftthreetimes}

\newcommand{\C}{\ensuremath{\mathbb C}}

\newcommand{\E}{\ensuremath{\mathbb E}}
\newcommand{\R}{\ensuremath{\mathbb R}}

\newcommand{\N}{\ensuremath{\mathbb N}}

\newcommand{\cR}{\mathcal R}
\newcommand{\cH}{\mathcal H}

\renewcommand{\l}{\lambda}
\renewcommand{\a}{\alpha}

\def\Hor{\mathop{\text{Hor}}}
\def\gr{\mathop{\text{gr}}}

\def\SO{\mathop{\text{SO}}}

\renewcommand\O{\mathop{\text{O}}}

\newcommand{\phil}{\varphi_{\l}}

\begin{document}
\makeatletter
\title{The dual horospherical Radon transform for polynomials}
\author{J. Hilgert\thanks{
Institut f\"ur Mathematik, Technische Universit\"at Clausthal,
38678 Clausthal-Zellerfeld, Germany.
E-mail: hilgert@math.tu-clausthal.de},
\/  A. Pasquale\thanks{%
Institut f\"ur Mathematik, Technische Universit\"at Clausthal,
38678 Clausthal-Zellerfeld, Germany.
E-mail: mapa@math.tu-clausthal.de}
\/ and E.B. Vinberg\thanks{%
Department of Algebra, Moscow State University, 119899 Moscow, Russia.
E-mail:vinberg@ebv.pvt.msu.su}
\thanks{The work of the third author was partially supported by CRDF
grant RM1-2088}}
\date{}
\maketitle
\makeatother

\section{Introduction}

A Radon transform is generally associated to a
double fibration

\bigskip
\begin{center}
\setlength{\unitlength}{4100sp}%
\begingroup\makeatletter\ifx\SetFigFont\undefined%
\gdef\SetFigFont#1#2#3#4#5{%
  \reset@font\fontsize{#1}{#2pt}%
  \fontfamily{#3}\fontseries{#4}\fontshape{#5}%
  \selectfont}%
\fi\endgroup%
\begin{picture}(1734,564)(79,242)
\thinlines
\put(800,780){$Z$}
\put(50,180){$X$}
\put(1500,180){$Y$,}
\put(300,600){$p$}
\put(1350,600){$q$}
\put(700,700){\vector(-4,-3){400}}
\put(1000,700){\vector(4,-3){400}}
\end{picture}
\end{center}
\bigskip
where one may assume without loss of generality that the maps
$p$ and $q$ are surjective and $Z$ is embedded into $X\times Y$
via $z\mapsto \big(p(z),q(z)\big)$. Let some measures be chosen
on $X,Y,Z$ and on the fibers of $p$ and $q$ so that
\begin{equation}
\label{(1)}
\int_X\left(\int_{p^{-1}(x)}f(u)du\right)dx
            =\int_Zf(z)\,dz=\int_Y\left(\int_{q^{-1}(y)}f(v)dv\right)dy.
\end{equation}
Then the Radon transform $\cR$ is the linear map
assigning to a function $\varphi$ on $X$
the function on $Y$ defined by
$$(\cR \varphi)(y)=\int_{q^{-1}(y)} (p^*\varphi)(v)\,dv,$$
where we have set $p^*\varphi:=\varphi \circ p$.
In a dual fashion, one defines a linear
transform $\cR^*$ from functions on $Y$
to functions on $X$ via
$$(\cR^*\psi)(x)=\int_{p^{-1}(x)} (q^*\psi)(u)\,du.$$
It is dual to $\cR$. Indeed, formally, \begin{align*} (\cR\varphi,\psi)
:&=\int_Y\left(\int_{q^{-1}(y)}(p^*\varphi)(v)\,dv\right)\psi(y)\,dy\\
&=\int_Y\left(\int_{q^{-1}(y)}(p^*\varphi)(q^*\psi)(v)\,dv\right)dy\\
&=\int_Z(p^*\varphi)(z)(q^*\psi)(z)\,dz \displaybreak[0] \\
&=\int_X \varphi(x) \left(\int_{p^{-1}(x)}(q^*\psi)(u)\,du\right)dx\\
&=:(\varphi,\cR^*\psi).
\end{align*}

In particular, if $X=G/K$ and $Y=G/H$ are homogeneous spaces of a Lie group
$G$, one can take $Z=G/(K\cap H)$, with $p$ and $q$ being $G$-equivariant maps
sending $e(K\cap H)$ to $eK$ and $eH$, respectively. Assume that there
exist $G$-invariant measures on $X$, $Y$ and $Z$. If such measures are fixed,
one can uniquely define measures on the fibers of $p$ and $q$ so that
the condition (\ref{(1)}) holds. (Here, speaking about  a measure on a
smooth manifold, we mean a measure defined by a differential form of top
degree.) In this setting we can consider the transforms $\cR$ and $\cR^*$,
and if $\dim X=\dim Y$, one can hope that they are invertible.  The basic
example is provided by the classical Radon transform, which acts on
functions on the Euclidean space $\E^n=(\R^n\sdir\SO(n))/\SO(n)$ and maps
them to functions on the space $H\E^n=(\R^n\sdir
\SO(n))/(\R^{n-1}\sdir\O(n-1))$ of hyperplanes in $\E^n$.

For a semisimple Riemannian symmetric space $X=G/K$ of noncompact type,
one can consider the horospherical  Radon transform as proposed by
I.M. Gelfand and M.I. Graev in \cite{GG59}. Namely, generalizing  the classical
notion of a horosphere in Lobachevsky space, one can define a horosphere
in $X$ as an orbit of a maximal unipotent subgroup of $G$. The group $G$
naturally acts on the set $\Hor X$ of all horospheres. This action
is transitive, so we can identify $\Hor X$ with some quotient space
$G/S$ (see \S \ref{section:horos} for details). It turns out that
$\dim X=\dim \Hor X$. Moreover, the groups $G$, $K$, $S$ and $K\cap S=M$
are unimodular, so there exist $G$-invariant measures on $X$, $\Hor X$
and $Z=G/M$. The Radon transform $\cR$ associated to the double fibration

\bigskip
\begin{center}
\setlength{\unitlength}{4100sp}%
\begingroup\makeatletter\ifx\SetFigFont\undefined%
\gdef\SetFigFont#1#2#3#4#5{%
  \reset@font\fontsize{#1}{#2pt}%
  \fontfamily{#3}\fontseries{#4}\fontshape{#5}%
  \selectfont}%
\fi\endgroup%
\begin{picture}(1734,564)(79,242)
\thinlines
\put(700,780){$G/M$}
\put(-200,150){$X=G/K$}
\put(1250,150){$G/S=\Hor X$}
\put(300,600){$p$}
\put(1450,600){$q$}
\put(600,700){\vector(-4,-3){400}}
\put(1200,700){\vector(4,-3){400}}
\end{picture}
\end{center}
\bigskip
is called the \emph{horospherical Radon transform}.

The space $Z=G/M$ can be interpreted as the set of pairs
$(x,\cH)\in X\times \Hor X$
with $x\in\cH$, so that $p$ and $q$ are just the natural projections. The
fiber $q^{-1}(\cH)$ with $\cH\in \Hor X$ is then identified with the
horosphere $\cH$, and the fiber $p^{-1}(x)$ with $x\in X$ is identified
with the submanifold $\Hor_x X\subseteq \Hor X$ of all horospheres passing
through $x$.  Note that, in contrast to the horospheres, all submanifolds
$\Hor_x X$ are compact, since $\Hor_x X$ is the orbit of the stabilizer of
$x$ in $G$, which is conjugate to $K$.

In this paper, we describe the \emph{dual horospherical Radon transform}
$\cR^*$ in terms of its action on
polynomial functions.
Here a differentiable function $\varphi$ on a homogeneous space $Y=G/H$
of a Lie group $G$ is called {\it polynomial}, if the linear span of the
functions $g\varphi$ with $g\in G$ is finite dimensional. The polynomial
functions constitute an algebra denoted by $\R[Y]$.

For $X=G/K$ as above, the algebra $\R[X]$ is finitely generated and $X$
is naturally identified with a connected component of the corresponding
affine real algebraic variety (the real spectrum of $\R[X]$). The natural
linear representation of $G$ in $\R[X]$ decomposes into a sum of mutually
non-isomorphic absolutely irreducible finite dimensional representations
whose highest weights $\lambda$ form a semigroup $\Lambda$.
Let $\R[X]_\lambda$ be the irreducible component of $\R[X]$ with highest weight
$\lambda$, so
\begin{equation}
\label{(2)}
\R[X]=\bigoplus_{\lambda\in\Lambda}\R[X]_\lambda.
\end{equation}
Denote by $\varphi_\lambda$ the highest weight function in $\R[X]_\lambda$
normalized by the condition
$$\varphi_\lambda(o)=1,$$
where $o=eK$ is the base point of $X$.
Then the subgroup $S$ is the intersection of the stabilizers of
all $\varphi_\lambda$'s. Its
unipotent radical $U$ is a maximal unipotent subgroup of $G$.

The algebra $\R[\Hor X]$ is also finitely generated. The manifold $\Hor X$
is naturally identified with a connected component of a quasi-affine algebraic
variety, which is a Zariski open subset in the real spectrum of $\R[\Hor X]$.
The natural linear representation of $G$ in $\R[\Hor X]$ is isomorphic to the
representation of $G$ in $\R[X]$.
Let $\R[\Hor X]_\lambda$ be the irreducible component of $\R[\Hor X]$ with
highest weight $\lambda$, so that
\begin{equation}
\label{(3)}
\R[\Hor X]=\bigoplus_{\lambda\in \Lambda} \R[\Hor X]_\lambda.
\end{equation}
Denote by $\psi_\lambda$ the highest weight function in $\R[\Hor X]_\lambda$
normalized by the condition
$$\psi_\lambda(sUo)=1,$$
where $s$ is the symmetry with respect to $o$ and
$sUo=(sUs^{-1})o$ is considered as a point of $\Hor X$.

The decomposition (\ref{(2)}) defines a filtration of the algebra
$\R[X]$ (see \S \ref{section:reps} for the precise definition).
Let $\gr\R[X]$ be the associated graded algebra. There is a canonical
$G$-equivariant algebra isomorphism
$$\Gamma\colon \gr\R[X]\to \R[\Hor X]$$
mapping each $\varphi_\lambda$ to $\psi_\lambda$. As a $G$-module,
$\gr \R[X]$ is canonically identified with $\R[X]$, so we can
view $\Gamma$ as a $G$-module isomorphism from $\R[X]$ to $\R[\Hor X]$.

The horospherical Radon transform $\cR$ is not defined for polynomial
functions on $X$ but its dual transform $\cR^*$ is defined for polynomial
functions on $\Hor X$, since it reduces to integrating along compact
submanifolds. Moreover, as follows from the definition of polynomial
functions, it maps polynomial functions on $\Hor X$ to polynomial
functions on $X$. Obviously, it is $G$-equivariant.
Thus, we have $G$-equivariant linear maps
$$\R[X]\overset{\Gamma}{\longrightarrow}\R[\Hor X]
       \overset{\cR^*}{\longrightarrow}\R[X].$$
Their composition $\cR^*\circ \Gamma$ is a $G$-equivariant linear operator
on $\R[X]$, so
$$(\cR^*\circ \Gamma)(\varphi)=
c_\lambda \varphi\quad \forall \varphi\in \R[X]_\l,$$
where the $c_\lambda$ are constants. To give a complete description of
$\cR^*$, it is therefore sufficient to find these constants. Our main result
is the following theorem.

\begin{Thm}
$c_\lambda=\mathbf{c}(\lambda+\rho)$, where $\mathbf{c}$ is the Harish-Chandra
$\mathbf{c}$-function and $\rho$ is the half-sum of the positive roots of $X$
(counted with multiplicities).
\end{Thm}

The Harish-Chandra $\mathbf{c}$-function
governs the asymptotic behavior of the zonal spherical functions on $X$.
A product formula for the
$\mathbf{c}$-function was found by S.G.\ Gindikin and F.I.\ Karpelevich
[GK62]: for a rank-one Riemannian symmetric space of the noncompact type,
the $\mathbf{c}$-function is a ratio of gamma functions
involving only the root multiplicities;
in the general case, it is the product of the $\mathbf{c}$-functions for
the rank-one
symmetric spaces defined by the indivisible roots of the space.
Thus, known the root structure of the symmetric space, the product formula
makes the $\mathbf{c}$-function, and hence our description of the
dual horospherical Radon transform, explicitly computable.

For convenience of the reader, we collect some crucial facts about the
$\mathbf{c}$-function in an appendix to this paper.

\bigskip

The following basic notation will be used in the paper without further
comments.
\begin{itemize}
\item Lie groups are denoted by capital Latin letters, and
      their Lie algebras by the corresponding small Gothic letters.
\item The dual space of a vector space $V$ is denoted by $V^*$.
\item The complexification of a real vector space $V$ is denoted by $V(\C)$.
\item The centralizer (resp. the normalizer) of a subgroup $H$ in a group
$G$ is denoted by $Z_G(H)$ (resp. $N_G(H)$).
\item The centralizer (resp. the normalizer) of a subalgebra $\frak h$ in
a Lie algebra $\frak g$ is denoted by $\frak z_{\frak g}(\frak h)$
(resp. by $\frak n_{\frak g}(\frak h)$).
\item If a group $G$ acts on a set $X$, we denote by $X^G$ the subset
of fixed points of $G$ in $X$.
\end{itemize}

%%%%%%%%%%%%%%%%%%%%%%%%%%%%%%%%%%%%%%%%%%%%

\section{Groups, spaces, and functions} \label{section:groups}

For any connected semisimple Lie group $G$ admitting a faithful (finite-dimensional)
linear representation, there is a connected complex algebraic group
defined over $\R$ such that $G$ is the connected component
of the group of its real points.
Among all such algebraic groups, there is a unique one such that all
the others are
its quotients. It is called the \emph{complex hull} of $G$ and is
denoted by $G(\C)$.
The group of real points of $G(\C)$ is denoted by $G(\R)$.
If $G(\C)$ is simply connected, then $G=G(\R)$.

The restrictions of polynomial functions on the algebraic group $G(\R)$
to $G$ are called {\it polynomial functions} on $G$. They are  precisely
those differentiable functions $\varphi$ for which the linear span
of the functions $g\varphi$, $g\in G$, is finite dimensional (see
e.g. \cite{CSM}, \S II.8).
(Here $G$ is supposed to act on itself by left multiplications.) The
polynomial functions on $G$ constitute an algebra which we denote by
$\R[G]$ and which is naturally isomorphic to $\R[G(\R)]$.

For any subgroup $H\subseteq G$, we denote by $H(\C)$ (resp. $H(\R)$)
its Zariski closure in $G(\C)$ (resp. $G(\R)$).
If $H$ is Zariski closed in $G$, i.e.
$H=H(\R)\cap G$, then $H$ is a subgroup of finite index in $H(\R)$; if $H$
is a semidirect product of a connected unipotent group and a compact group,
then $H=H(\R)$.

For a homogeneous space $Y=G/H$ with $H$ Zariski closed in $G$,
set $Y(\C)=G(\C)/H(\C)$. This is an algebraic variety defined over $\R$,
and $Y$ is naturally identified with a connected component of the variety
$Y(\R)$ of real points of $Y(\C)$. We call $Y(\C)$ the {\it complex hull}
of $Y$.

If $H$ is reductive, then also $H(\C)$ is reductive and the variety
$Y(\C)$ is affine, the algebra $\C[Y(\C)]$ being naturally
isomorphic to the algebra $\C[G(\C)]^{H(\C)}$ of $H(\C)$-right-invariant
polynomial functions on $G(\C)$
(see e.g. [VP89], Section 4.7 and Theorem 4.10). Correspondingly,
the algebra $\R[Y(\R)]$ is naturally isomorphic to $\R[G(\R)]^{H(\R)}$.

In general, the functions on $Y(\R)$ arising from  $H(\R)$-right-invariant
polynomial functions on $G(\R)$ are called polynomial functions on $Y(\R)$,
and their restrictions to $Y$ are called {\it polynomial functions} on $Y$.
They are precisely those differentiable functions $\varphi$ for which the
linear span of the functions $g\varphi$, $g\in G$, is finite dimensional.
They form an algebra which we denote by $\R[Y]$.

In the following we consider a semisimple Riemannian symmetric space $X=G/K$
of noncompact type. This means that $G$ is a connected semisimple Lie group
without compact factors and $K$ is a maximal compact subgroup of $G$.
We do not assume that the center of $G$ is trivial, so the action of $G$ on $X$
may be non-effective. We do, however, require that $G$ has a faithful linear
representation. According to the above, the space $X$ is then
a connected component of the affine algebraic variety $X(\R)$.

%%%%%%%%%%%%%%%%%%%%%%%%%%%%%%%%%%%%%%%%%%

\section{Subgroups and subalgebras} \label{section:subgroups}

We recall some facts about the structure of Riemannian symmetric
spaces of noncompact type (see \cite{He78} for details).
Let $X=G/K$ be as above and $\theta$ be the Cartan involution
of $G$ with respect to $K$, so $K=G^\theta$. Let $\frak a$ be a Cartan
subalgebra for $X$, i.e. a maximal abelian subalgebra in the
$(-1)$-eigenspace of $d\theta$.
 Its dimension $r$ is called the
\emph{rank} of $X$.
Under any representation of $G$, the elements of $\frak a$ are simultaneously
diagonalizable.
The group $A=\exp \frak a$ is  a
maximal connected abelian subgroup of $G$ such that $\theta(a)=a^{-1}$
for all $a \in A$. It is isomorphic to $(\R^*_+)^r$. Its Zariski closure
$A(\R)$ in $G(\R)$ is a split algebraic torus which is isomorphic to
$(\R^*)^r$.
Let $\scr X(A)$ denote the (additively written) group of real characters of
the torus $A(\R)$. It is a free abelian group of rank $r$.
We identify each character $\chi$ with its differential $d\chi \in \frak a^*$.

%For subgroups $B$ and $C$ of $G$, we write
%$Z_C(B)$ (resp. $N_C(B)$) for the centralizer (resp. normalizer) of $B$ in $C$.
%If $\frak b$ and $\frak c$ are subsets of $\frak g$,
%we also set
%$\frak z_{\frak c}(\frak b)=\{Y \in \frak c: [Y,Z]=0 \ \forall Z \in \frak b\}$
%for the centralizer of $\frak b$ in $\frak c$.

The root decomposition of $\frak g$ with respect to $A$
(or with respect to $A(\R)$, which is the same)
 is of the form
$$
\frak g=\frak g_0 +\sum_{\a \in \Delta} \frak g_\a,
$$
where $\frak g_0=\frak z_{\frak g}(\frak a)$.
%(For subalgebras $\frak b$ and
%$\frak c$ of $\frak g$  we use the notation
%$\frak z_{\frak c}(\frak b):=\{Y \in \frak c: [Y,\frak b]=0\}$
%for the centralizer of $\frak b$ in $\frak c$.)
If $\frak m:=\frak z_{\frak k}(\frak a)$, then $\frak g_0=\frak m +\frak a$.

The set $\Delta \subset \scr X(A)$ is the \emph{root system} of $X$
(or the \emph{restricted root system} of $G$) with respect to $A$
and $\frak g_\a$ is the \emph{root subspace}
corresponding to $\alpha$. The dimension of $\frak g_\a$ is called the
\emph{multiplicity} of the root $\a$ and is denoted by $m_\a$.
By the identification of a character with its differential, we will consider
$\Delta$ as a subset of $\frak a^*$.
Choose a system $\Delta^+$ of \emph{positive roots} in $\Delta$. Let $\Pi=\{\a_1,\dots,\a_r\}
\subset \Delta^+$ be the corresponding system of \emph{simple roots}. Then
$$C=\{x \in \frak a: \a_i(x) \ge 0 \ \text{for $i=1,\dots,r$}\}$$
is called the \emph{Weyl chamber} with respect to $\Delta^+$.
The subspace
$$\frak u=\sum_{\a \in \Delta^+} \frak g_\a$$
is a \emph{maximal unipotent subalgebra} of $\frak g$.

Set
$$G_0=Z_G(A), \quad M=Z_K(A).$$
%(If  $B$ and $C$ are subgroups of $G$, we write
%$Z_C(B)$ for the centralizer of $B$ in $C$ and
%$N_C(B)$ for the normalizer of $B$ in $C$.)
Then $G_0=M\times A$ and the Lie algebras of $G_0$ and $M$ are
$\frak g_0$ and $\frak m $, respectively.
Clearly, $G_0$ is Zariski closed in $G$.

The group $U=\exp \frak u$ is a
\emph{maximal unipotent subgroup} of $G$. It is normalized by $A$ and the map
$$U \times A \times K \rightarrow G, \quad (u,a,k) \mapsto uak,$$
is a diffeomorphism. The decomposition $G=UAK$ (or $G=KAU$) is called the \emph{Iwasawa
decomposition} of $G$.
Since every root subspace is $G_0$-invariant, $G_0$ normalizes $U$, so
 $$P:=UG_0=U \leftthreetimes G_0$$
is a subgroup of $G$. Moreover, $P=N_G(U)$
(see e.g. \cite{Wa72}, Proposition 1.2.3.4), so $P$ is Zariski
closed in $G$.

We say that a Zariski closed subgroup of $G$ is \emph{parabolic}, if its Zariski closure in
$G(\C)$ is a parabolic subgroup of $G(\C)$.
Then $P$ is a \emph{minimal parabolic subgroup} of $G$.
The subgroup
$$S:=UM=U \leftthreetimes M$$
of $P$ is normal in $P$ and $P/S$ is isomorphic to $A$.
It follows from the Iwasawa decomposition that $K\cap S=M$.

%%%%%%%%%%%%%%%%%%%%%%%%%%%%%%%%%%%%%%%%%%%%%%%
\section{Representations} \label{section:reps}

For later use we collect some well-known facts about
finite-dimensional representations of $G$.

The natural linear representation of $G$ in $\R[X]$ decomposes into a sum of
mutually non-isomorphic absolutely irreducible finite-dimensional
representations called (finite-dimensional) \emph{spherical
representations} (see e.g. \cite{GW98}, chap.\ 12).

%If $G$ acts on a vector space $V$ and $H$ is a subgroup of $G$,
%we denote by $V^H$ the
%space of vectors in $V$ fixed by $H$.

\begin{Thm} {\rm (see \cite{He84}, \S\ V.4)} An irreducible finite dimensional
representation of
$G$ on a real vector space $V$ is spherical if and only if the following
equivalent conditions hold:
\begin{enumerate}
\item[{\rm (1)}] $V^K\not=\{0\}$.
\item[{\rm (2)}] $V^S\not=\{0\}$.
\end{enumerate}
If these conditions hold, then $\dim V^K=\dim V^S=1$ and the subspace
$V^S$ is invariant under $P$.
\end{Thm}

For a spherical representation, the group $P$ acts on $V^S$ via
multiplication by some character of $P$ vanishing on $S$. The restriction
of this character to $A$ is called the {\it highest weight} of the
representation. A spherical representation is uniquely determined
by its highest weight.
The highest weights of all irreducible spherical representations constitute a
subsemigroup $\Lambda \subset \scr X(A)$.

An explicit description of
$\Lambda$ as a subset of $\frak a^*$ can be given as follows.
Let $\Delta_\star$ denote the set of roots $\a \in \Delta$ such that
$2\a \notin \Delta$. Then $\Delta_\star$ is a root system in $\frak a^*$,
and a system of simple roots
corresponding to $\Delta_\star^+:=\Delta^+ \cap \Delta_*$ can be obtained
from the system $\Pi=\{\a_1,\dots, \a_r\}$ of simple roots in $\Delta^+$
by setting for $j=1,\dots, r$
\begin{equation*}
\beta_j:= \begin{cases} \a_j &\text{if $2\a_j \notin \Delta^+$},\\
                                   2\a_j &\text{if $2\a_j \in \Delta^+$}.
              \end{cases}
\end{equation*}

Let $\omega_1,\dots, \omega_r \in \frak a^*$ be defined by
\begin{equation} \label{eq:mu}
\frac{(\omega_j,\beta_i)}{(\beta_i,\beta_i)}=\delta_{ij},
\end{equation}
where $( \cdot , \cdot )$ denotes the scalar product in $\frak a^*$
induced by an invariant scalar product in $\frak g$.

\begin{Prop}  \label{prop:weights}
{\rm (see \cite{He94}, Proposition 4.23)}
The semigroup $\Lambda$ is freely generated by $\omega_1,\dots, \omega_r$.
\end{Prop} \qed

The weights $\omega_1,\dots, \omega_r$ are called the \emph{restricted
weights} of $X$.

For $\lambda \in \Lambda$, let $\R[X]_\lambda$ denote
the irreducible component of $\R[X]$
with highest weight $\lambda$. Then we have
$$
\R[X]=\bigoplus_{\lambda \in \Lambda} \R[X]_\lambda.
$$
We denote by $\phil^S$ the highest weight function of $\R[X]_\lambda$
normalized by the condition
$$\phil^S(o)=1.$$
(Since $Po=PKo=Go=X$,
we have $\phil^S(o)\neq 0$.)
Obviously,
$$
\varphi_\lambda^S \varphi_\mu^S=\varphi_{\lambda+\mu}^S.
$$

In general, the multiplication in $\R[X]$ has the property
$$
\R[X]_\lambda \R[X]_\mu \subset \bigoplus_{\nu \leq \lambda+\mu} \R[X]_\nu,$$
where ``$\leq$'' is the ordering in the group ${\scr X}(A)$
defined by the subsemigroup
generated by the simple roots.
In other words, the subspaces
$$
\R[X]_{\leq \lambda} =\bigoplus_{\mu \leq \lambda} \R[X]_\mu
$$
constitute a ${\scr X}(A)$-filtration of the algebra $\R[X]$ with respect to the ordering ``$\leq$''.

The functions $\varphi \in \R[X]_\lambda$ vanishing at $o$ constitute a $K$-invariant subspace of codimension $1$.
The $K$-invariant complement of it is a $1$-dimensional subspace,
on which $K$ acts trivially.
Let $\varphi_\lambda^K$ denote the function of this subspace normalized
by the condition $\phil^K(o)=1$. It is called the \emph{zonal spherical function of
weight $\lambda$}.

\begin{Lemma} {\rm (see e.g. \cite{He84}, p. 537)}\label{invscal}
For any finite dimensional irreducible representation of $G$ on a
real vector space
$V$ there is a positive definite scalar product
$(\cdot\mid\cdot)$ on $V$ such that
$$(gx\mid\theta(g)y)=(x\mid y) \quad \text{for all $g\in G$ and $x,y\in V$.}$$
This scalar product is unique up to a scalar multiple.
\end{Lemma}

The scalar product given by Lemma \ref{invscal} is called $G$-skew-invariant.
Note that it is $K$-invariant.

With respect to  the $G$-skew-invariant scalar product on $\R[X]_\lambda$
the zonal spherical function $\varphi_\lambda^K$ is orthogonal  to the subspace of
functions vanishing at $o\in X$. Let $\alpha_\lambda$ denote the angle between
$\varphi_\lambda^K$ and $\varphi_\lambda^S$. Then the projection of $\varphi_\lambda^K$ to
$\R[X]_\lambda^S=\R\varphi_\lambda^S$ is equal to $(\cos^2\alpha_\lambda)\varphi_\lambda^S$
(see Figure 1). In particular, we see that $\varphi_\lambda^K$ and $\varphi_\lambda^S$
are not orthogonal.
\begin{center}
\epsfxsize5cm
\epsfbox{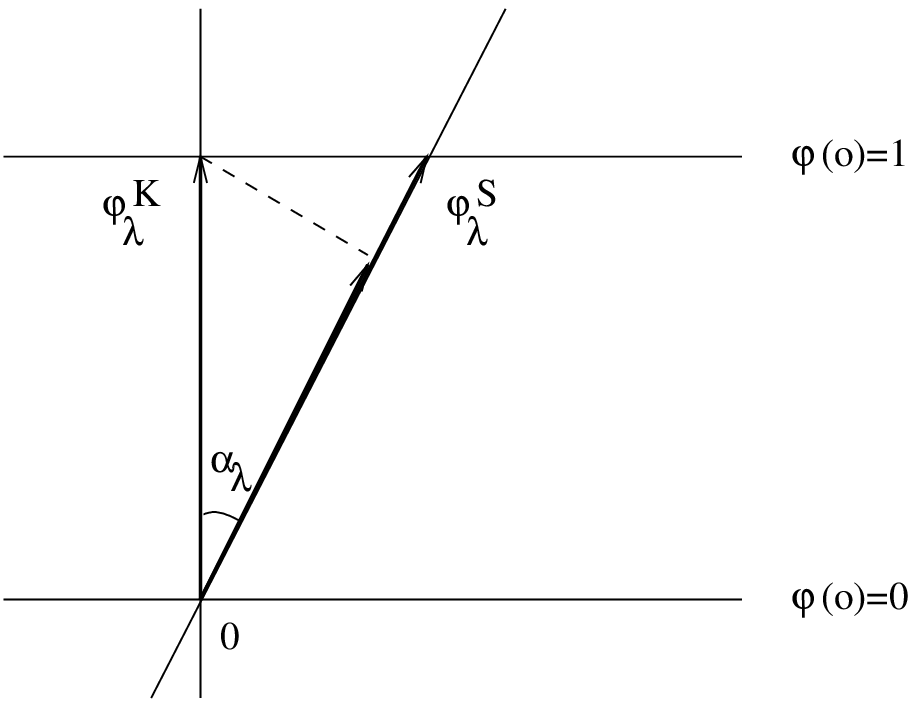}

Figure 1
\end{center}
The weight decomposition of $\varphi_\lambda^K$ is of the form
$$\varphi_\lambda^K=(\cos^2\alpha_\lambda)\varphi_\lambda^S
  +\sum_{\mu<\lambda}\varphi_{\lambda,\mu},$$
for some weight vector $\varphi_{\lambda,\mu}$ of weight $\mu$ in
$\R[X]_\lambda$.
This gives the asymptotic behavior of $\varphi_\lambda^K$ on
$(\exp(-C^o))o$, where $C^o$ is the interior of the Weyl chamber $C$ in $\frak a$.
More precisely, for $\xi\in C^o$ we have
$$\varphi_\lambda^K((\exp(-t\xi))o)=((\exp t\xi)\varphi_\lambda^K)(o)
  \underset{t\to+\infty}{\sim}
  (\cos^2\alpha_\lambda)e^{t\lambda(\xi)}.$$
But it is known (see \cite{He84}, \S~IV.6)
that the same asymptotics is described in
terms of the Harish-Chandra $\mathbf{c}$-function:
$$\varphi_\lambda^K((\exp(-t\xi))o)\underset{t\to+\infty}{\sim}
  \mathbf{c}(\lambda+\rho)e^{t\lambda(\xi)},$$
where $\rho=\frac{1}{2}\sum_{\alpha\in \Delta+}m_\alpha \alpha$
is the half-sum of positive roots.
This shows that
$$   \cos^2\alpha_\lambda=\mathbf{c}(\lambda+\rho).$$

%%%%%%%%%%%%%%%%%%%%%%%%%%%%%%%%%%%%%%%%%%%%
\section{Horospheres} \label{section:horos}

\begin{Def}
A \emph{horosphere} in $X$ is an orbit of a maximal unipotent subgroup of $G$.
\end{Def}

Since all maximal unipotent subgroups are conjugate to $U$, any horosphere
is of the form
$gUx$  ($g \in G$, $x \in X$).
Moreover, since $X=Po$ and $P$ normalizes $U$, any horosphere can be
represented in the form $gUo$ ($g \in G$). In other words,
the set $\Hor X$ of all horospheres is a homogeneous space of $G$.
The following lemma shows that the $G$-set $\Hor X$ is identified with $G/S$
if we take the horosphere $Uo$ as the base point for $\Hor X$.

\begin{Lemma}{\rm (see also \cite{He94}, Theorem 1.1, p.~77)}
The stabilizer of the horosphere $Uo$ is the algebraic subgroup
$$S=UM =U \leftthreetimes M.$$
\end{Lemma}
\pf
Obviously, $S$ stabilizes $Uo$.
Hence the stabilizer of $Uo$ can be written as $\widetilde S=U \widetilde
M$, where $$M \subset \widetilde M \subset K.$$ Since $U$ is a maximal
unipotent subgroup in $G$ (and hence in $\widetilde S$), it contains the
unipotent radical $\widetilde U$ of $\widetilde S$.  The reductive group
$\widetilde S/\widetilde U$ can be decomposed as $$\widetilde S/\widetilde
U=(U/\widetilde U)\widetilde M,$$ so the manifold $(\widetilde
S/\widetilde U)/(U/\widetilde U)$ is compact.  But then the Iwasawa
decomposition for $\widetilde S/\widetilde U$ shows that the real rank of
$\widetilde S/\widetilde U$ equals $0$, that is $\widetilde S/\widetilde
U$ is compact. This is possible only if $U=\widetilde U$. Hence, $$
\widetilde S \subset N(U)=P=S \leftthreetimes A.$$ It follows from the
Iwasawa decomposition $G=AUK$ that $\widetilde S \cap A= \{e\}$. Thus
$\widetilde S=S$.  \qed

As follows from \cite{VP72}, the $G$-module structure of $\R[\Hor X]$
is exactly  the same as that of $\R[X]$, but in contrast to the case of
$\R[X]$, the decomposition of $\R[\Hor X]$ into the sum of irreducible
components is a $\scr X(A)$-grading.

Let $\R[\Hor X]_\lambda$ be the irreducible component of $\R[\Hor X]$
with highest weight $\lambda$, so
$$\R[\Hor X]=\bigoplus_{\lambda\in \Lambda}\R[\Hor X]_\lambda.$$

Denote by $\psi_\lambda^S$ and $\psi_\lambda^K$ the highest weight
function and the $K$-invariant function in $\R[\Hor X]_\lambda$ normalized by
$$\psi_\lambda^S(sUo)=\psi_\lambda^K(sUo)=1.$$
To see that this is possible, note that the horosphere $sUo$ is stabilized by
$sSs^{-1}=\theta(S)$. Hence the subspace $V_0$ of functions in $\R[\Hor X]_\lambda$ vanishing
at $sUo$ is $\theta(S)$-invariant.
Its orthogonal complement is therefore $S$-invariant
and must coincide with $\R[\Hor X]_\lambda^S$. This implies
$\psi_\lambda^S(sUo)\not=0$, so we can normalize $\psi_\lambda^S$ as asserted.
Since $\R[\Hor X]_\lambda^K$ and $\R[\Hor X]_\lambda^K$ are not orthogonal in
$\R[\Hor X]_\lambda\cong \R[X]_\lambda$, we have
 $\R[\Hor X]_\lambda^K\cap V_0=\{0\}$
and we can
normalize also $\psi_\lambda^K$ as asserted.
Notice that the horospheres passing through $o$ form a single $K$-orbit and
therefore $\psi_\lambda^K$ takes the value $1$ on each of them.

Consider the ${\mathbf X}(A)$-graded algebra $\gr \R[X]$ associated with the
${\mathbf X}(A)$-filtration of $\R[X]$ defined in \S \ref{section:reps}.
As a $G$-module, $\gr\R[X]$ can be identified with $\R[X]$,
but when we multiply elements
$\varphi_\lambda\in \R[X]_\lambda$ and $\varphi_\mu\in \R[X]_\mu$
 in $\gr\R[X]$, only the highest term in
their product in $\R[X]$ survives.
Moreover, there is a unique $G$-equivariant linear isomorphism
$$\Gamma\colon \R[X]=\gr\R[X]\to \R[\Hor X]$$
such that $\Gamma(\varphi_\lambda^S)=\psi_\lambda^S$.

\begin{Prop} The map $\Gamma$ is an isomorphism of the algebra $\gr \R[X]$ onto the
algebra $\R[\Hor X]$.
\end{Prop}

\pf For any semisimple complex algebraic group, the tensor product
of the irreducible representations with highest weights $\lambda$ and
$\mu$ contains a unique irreducible component with highest weight $\lambda+\mu$.
It follows that, if we identify the irreducible components of $\R[X]$
with the corresponding irreducible components of $\R[\Hor X]$ via $\Gamma$,
the product of functions
$\varphi_\lambda\in \R[X]_\lambda$ and $\varphi_\mu\in \R[X]_\mu$
in $\gr\R[X]$ differs from
their product in $\R[\Hor X]$ only by some factor $a_{\lambda\mu}$
depending only on $\lambda$
and $\mu$. Taking $\varphi_\lambda=\varphi_\lambda^S$ and
$\varphi_\mu=\varphi_\mu^S$, we conclude that $a_{\lambda\mu}=1$.
\qed

\begin{Rem}\label{Gammanat}
The definition of $\Gamma$ makes use of the choice of a base point $o$ in
$X$ and a maximal unipotent subgroup $U$ of $G$, but it is easy to see
that all such pairs $(o,U)$ are $G$-equivalent. It follows that
$\Gamma$ is in fact canonically defined.
\end{Rem}

%%%%%%%%%%%%%%%%%%%%%%%%%%%%%%%%%%%%%%%%%
\section{Proof of the Main Theorem} \label{section:intgeom}

Consider the double fibration
\begin{center}
\setlength{\unitlength}{4100sp}%
\begingroup\makeatletter\ifx\SetFigFont\undefined%
\gdef\SetFigFont#1#2#3#4#5{%
  \reset@font\fontsize{#1}{#2pt}%
  \fontfamily{#3}\fontseries{#4}\fontshape{#5}%
  \selectfont}%
\fi\endgroup%
\begin{picture}(1734,564)(79,242)
\thinlines
\put(700,780){$G/M$}
\put(-200,150){$X=G/K$}
\put(1250,150){$G/S=\Hor X$}
\put(300,600){$p$}
\put(1450,600){$q$}
\put(600,700){\vector(-4,-3){400}}
\put(1200,700){\vector(4,-3){400}}
\end{picture}
\end{center}
\bigskip

\noindent
Since all the involved groups are unimodular, there are invariant measures on
the homogeneous spaces $X$, $\Hor X$, $G/M$ and on the fibers of $p$ and $q$,
which are the images under the action of $G$ of $K/M$ and $S/M$, respectively.
Let us normalize these measures so that:
\begin{enumerate}
\renewcommand{\labelenumi}{(\theenumi)}
\item
the volume of $K/M$ is $1$;
\item
the measure on $G/M$ is the product of the measures on $K/M$ and $X$;
\item
the measure on $G/M$ is the product of the measures on $S/M$ and
$\Hor X$.
\end{enumerate}
(This leaves two free parameters).

Consider the dual horospherical Radon transform
$$
\mathcal R^*: \R[\Hor X] \rightarrow \R[X].
$$
Combining it  with the map $\Gamma$ defined in \S \ref{section:horos},
we obtain a $G$-equivariant linear isomorphism
$$\mathcal R^*\circ \Gamma\colon \R[X]\to \R[X].$$
In view of absolute irreducibility, Schur's lemma shows that
$\mathcal R^*\circ \Gamma$ acts on each $\R[X]_\lambda$ by scalar
multiplication. The scalars are given by the following theorem:

\begin{Thm} \label{thm:bisc}
For $\varphi\in \R[X]_\lambda$,
$$(\mathcal R^* \circ \Gamma)(\varphi)=\mathbf c(\l+\rho)\varphi,$$
where $\mathbf c$ is the Harish-Chandra $\mathbf c$-function.
\end{Thm}
\pf
We test the map at the zonal spherical function $\varphi_\lambda^K\in \R[X]_\lambda$.
The map $\Gamma$ takes it to $c_\lambda\psi_\lambda^K$ for some $c_\lambda\in \R$.
Since the function $\psi_\lambda^K$ has value $1$ on the horospheres passing
through $o$, the
map $\mathcal R^*$ takes it to $\varphi_\lambda^K$. Thus we have
$$(\mathcal R^*\circ \Gamma)(\varphi_\lambda^K)=c_\lambda\varphi_\lambda^K.$$
Identifying $\R[X]_\lambda$ and $\R[\Hor X]_\lambda$ via $\Gamma$, we now  find
$c_\lambda=\cos^2\alpha_\lambda$ (see Figure 2), and this proves the claim.

\begin{center}
\epsfxsize5cm
\epsfbox{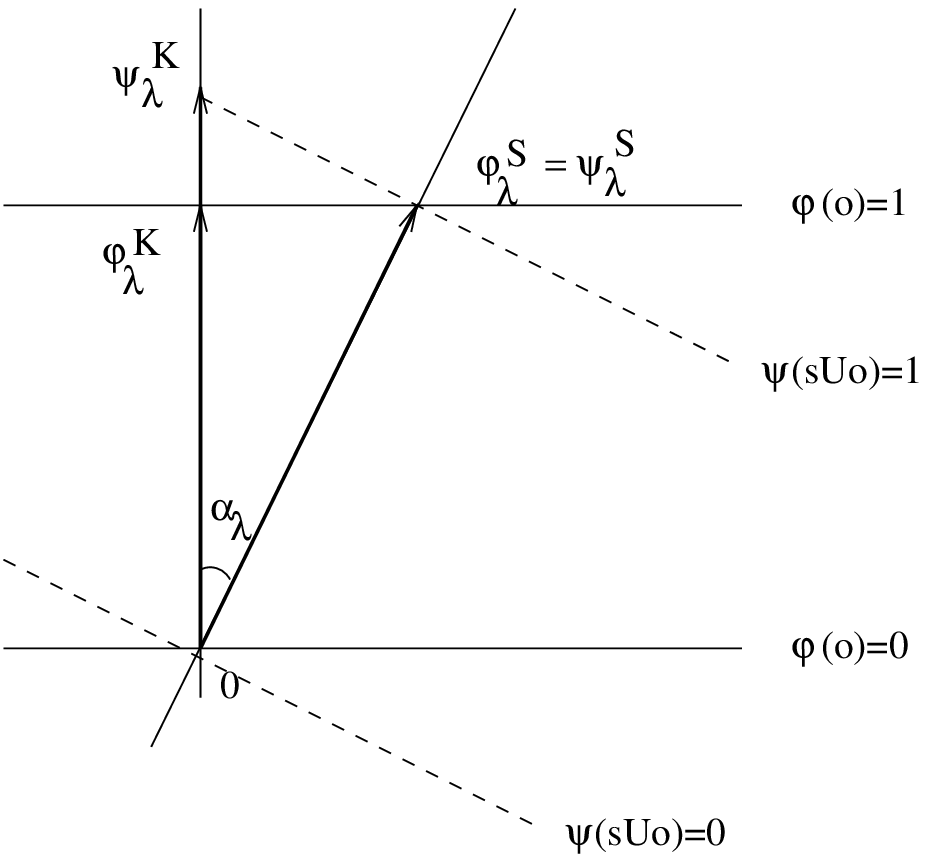}

Figure 2
\end{center}
\qed

%%%%%%%%%%%%%%%%%%%%%%%%%%%%%%%%%%%%%%%%%%%%%%%
\appendix

\section{Appendix: The $\mathbf c$-function} \label{section:c}

Because of the Iwasawa decomposition $G=KAU$,
every $g \in G$ can be written as
$g=k \exp H(g) u$ for a uniquely determined $H(g) \in \frak a$.
Let $\overline{U}:=\theta(U)$, and let $d\overline{u}$ denote the
invariant measure on $\overline{U}$ normalized by the condition
$$\int_{\overline{U}} e^{-2\rho(H(\overline{u}))}\; d\overline{u}=1.$$

The $\mathbf c$-function has been defined by Harish-Chandra as the
integral $$ \mathbf c(\l):=\int_{\overline{U}}
e^{-(\l+\rho)(H(\overline{u}))}\; d\overline{u}, $$ which absolutely
converges for all $\l \in \frak a(\C)^*$ satisfying ${\rm Re} (\l,\a)>0$
for all $\a \in \Delta^+$.  The computation of the integral gives the
so-called Gindikin--Karpelevich product formula \cite{GK62} (see also
\cite{He84}, Section IV.6.4, or \cite{GV}, p. 179):
\begin{equation} \label{eq:c}
\mathbf c(\l)= \kappa \prod_{\a \in \Delta^{++}}
\frac
{2^{-\l_\a} \; \Gamma\left(\l_\a \right)}
{\Gamma\left( \displaystyle{
\frac{\l_\a}{2} +\frac{m_\a}{4}+\frac{1}{2}}\right)
\Gamma\left( \displaystyle{
\frac{\l_\a}{2} +\frac{m_\a}{4}+\frac{m_{2\a}}{2}}\right)},
\end{equation}
where $\Delta^{++}$ denotes the set of indivisible roots in $\Delta^+$,
$\l_\a:=\frac{(\l,\a)}{(\a,\a)}$, and
the constant $\kappa$ is chosen so that $\mathbf c(\rho)=1$.
This formula provides the explicit meromorphic continuation of
$\mathbf c$ to the entire $\frak a(\C)^*$.

Formula (\ref{eq:c}) simplifies in the case of a reduced root system (i.e.
when $\Delta^{++}=\Delta^+$), since the duplication formula
\begin{equation} \label{eq:duplgamma}
\Gamma(2z)=2^{2z-1} \sqrt{\pi}\, \Gamma(z) \Gamma(z+1/2)
\end{equation}
for the gamma function yields
\begin{equation} \label{eq:creduced}
\mathbf c(\l)= \kappa \prod_{\a \in \Delta^{+}}
 \frac{\Gamma(\l_\a)}{\Gamma(\l_\a+m_\a/2)},
\end{equation}
with
\begin{equation*}
\kappa=\prod_{\a \in \Delta^{+}} \frac{\Gamma(\rho_\a+m_\a/2)}
{\Gamma(\rho_\a)}.
\end{equation*}

If, moreover, all the multiplicities $m_\a$ are even (which is
equivalent to the property that all Cartan subalgebras of $\frak g$
are conjugate), then
the functional equation $z \Gamma(z)=\Gamma(z+1)$ implies
\begin{equation*} \mathbf c(\l)= \prod_{\a \in \Delta^{+}}
\frac{\rho_\a(\rho_\a+1)\dots(\rho_\a+m_\a/2-1)}
{\l_\a(\l_\a+1)\dots(\l_\a+m_\a/2-1)}.
\end{equation*}

Finally, suppose that the group $G$ admits a complex structure.
In this case the root system is reduced
and $m_\a=2$ for every root $\a$, and (\ref{eq:creduced}) reduces to
\begin{equation*}
\mathbf c(\l)=
\prod_{\a \in \Delta^{+}} \frac{\rho_\a}{\l_\a}.
\end{equation*}

\begin{Ex}\label{ex:hyper}
For $n$-dimensional Lobachevsky
space, there is only one positive root $\a$, with $m_\a=n-1$, so
$\rho=(n-1)\a/2$. Let $\l=l\a$. Then formula (\ref{eq:creduced})
gives

$$
\mathbf c(\l)=\frac{\Gamma(n-1)\Gamma(l)}{\Gamma(\frac {n-1}{2})
\Gamma(l+\frac{n-1}{2})}.
$$

The semigroup $\Lambda$ is generated by $\a$. For $\l=l\a\, (l\in \N)$,
we obtain $$c_\l=\mathbf c(\l+\rho)=
\frac{\Gamma(n-1)\Gamma(l+\frac{n-1}{2})}{\Gamma(\frac {n-1}{2})
\Gamma(l+n-1)}=
\begin{cases}
\displaystyle{\frac{(n+l-1)(n+l)\cdots (n+2l-2)}
{2^{2l}\;\frac{n}2(\frac{n}2+1)\cdots (\frac{n}2+l-1)}}
    , &\text{$n$ even},\\
\noalign{\smallskip}
\displaystyle{\frac{(\frac{n-1}2+1)(\frac{n-1}2+2)\cdots (\frac{n-1}2+l-1)}
     {2n(n+1)\cdots (n+l-2)}}, &\text{$n$ odd}.
\end{cases}
$$
\end{Ex}

 %%%%%%%%%%%%%%%%%%%%%%%%%%%%%%%%%%%%%%%%%%%%%%%%%
%%%%%%%%%%%%%%%%%   Bibliography %%%%%%%%%%%%%%%%%%%%%%%%

\end{document}